\documentclass[article,12pt,3p]{elsarticle}
\biboptions{numbers,sort&compress} 
\usepackage[
singlelinecheck=false 
]{caption}
\usepackage{enumitem}
\usepackage{bm}
\usepackage{setspace}
\usepackage{makecell}
\usepackage{multirow}
\usepackage[dvipsnames,usenames]{xcolor}
\usepackage{amsmath}
\usepackage{amssymb}
\usepackage{lipsum}
\usepackage{comment}
\usepackage{xcolor}
\usepackage{subfig}
\usepackage{subcaption}

\journal{\mbox{ \ }}

\begin{document}
	\makeatletter
	\def\ps@pprintTitle{%
		\let\@oddhead\@empty
		\let\@evenhead\@empty
		\def\@oddfoot{\reset@font\hfil}%
		\let\@evenfoot\@oddfoot}
	\makeatother
	
	\renewcommand{\arraystretch}{1.3}

	\begin{frontmatter}
		
		\title{Constrained volume-difference site percolation model on the square lattice}
			
		\author[label1]{Charles S. do Amaral\corref{cor1}}
		\address[label1]{Departamento de Matem\'atica - Centro Federal de Educa\c c\~ao Tecnol\'ogica de Minas Gerais, \linebreak Av. Amazonas 7675, Belo Horizonte, MG, Brazil, Postal Code: 30.150-000 }
		\ead{charlesmat@cefetmg.br}
		\cortext[cor1]{Corresponding author.}
  		
		\date{\today}
		
		\begin{abstract}
   We study a percolation model with restrictions on the opening of sites on the square lattice. In this model, each site $s \in \mathbb{Z}^{2}$ starts closed and an attempt to open it occurs at time $t=t_s$, where $(t_s)_{s \in \mathbb{Z}^2}$ is a sequence of independent random variables uniformly distributed on the interval $[0,1]$. The site will open if the volume difference between the two largest clusters adjacent to it is greater than or equal to a constant $r$ or if it has at most one adjacent cluster. Through numerical analysis, we determine the critical threshold $t_c(r)$ for various values of $r$, verifying that $t_c(r)$ is non-decreasing in $r$ and that there exists a critical value $r_c=5$ beyond which percolation does not occur. Additionally, we find that the correlation length exponent of this model is equal to that of the ordinary percolation model. For $t = 1$ and $1 \leq r \leq 9$, we estimate the averages of the density of open sites, the number of distinct cluster volumes, and the volume of the largest cluster.
		\end{abstract}
		
	\end{frontmatter}

\section{Introduction}

	Percolation models have numerous applications in various scientific fields such as ecology, chemistry, social sciences, and biology \cite{konai,davis,ziff_app,chen,xie,huth,sole,bunde,webb}. The first percolation model was introduced by Broadbent and Hammersley in 1957 \cite{broadbent} with the idea of modeling the flow of a deterministic fluid through a random environment. In this model, we assign probabilities $p$ and $1-p$, independently, so that each site of an infinite and connected graph is \textit{open} and \textit{closed}, respectively. If a site is open, then fluid can pass through it; otherwise, the fluid passage is blocked. Each set of connected open sites is called \textit{cluster}. This model presents a geometric phase transition characterized by the emergence of an infinite cluster when $p$ exceeds a specific constant $p_c$ (\textit{percolation threshold} or \textit{critical point}). In this case, we say that \textit{percolation} has occurred. We will denote by $\mathcal{C}_s$ the cluster that contains the site $s$ and by $|\mathcal{C}_s|$ its \textit{volume} (number of sites that belong to the cluster $\mathcal{C}_s$).
 
    An analogous model can be formulated by assigning probabilities for each bond, rather than for each site, to be open or closed. These models are called \textit{Ordinary Site Percolation model} (OSPM) and \textit{Ordinary Bond Percolation model} (OBPM), respectively.
 
    Another way to define the OSPM is by assigning a random number $t_s$, chosen uniformly in $[0,1]$, to each site $s$. We then consider a time variable $t \in [0,1]$ that varies continuously. When $t=0$, all sites are closed, and each site $s$ opens at $t=t_s$. In this definition, we have that percolation occurs when $t>t_c=p_c$ (\textit{critical time}). Note that there is no spatial restriction on the opening of each site. The OBPM can also be defined this way.
	
	There are variations of this model that consider restrictions on the opening of sites or bonds. Motivated by the sol-gel transitions, Aldous proposes the \textit{Frozen Percolation model} \cite{aldous}. In this model, we sequentially open the bonds of a graph with the restriction that no cluster exceeds a critical volume. When this occurs, the bonds around the cluster \textit{freeze} and cannot open later.
 
    Gaunt proposed a model called \textit{Percolation with Restricted-Valence} \cite{gaunt0} in which the opening of a site can only occur if the number of adjacent open sites is less than a predetermined constant. The \textit{Constrained-degree Percolation Model} \cite{bnb} is analogous to this model, the difference is that it considers the opening of the bonds instead of the sites. These models are related to the study of dimers and polymers \cite{gaunt1, lang, furlan, soteros, wilkinson}.

    Other mathematical studies on variations of percolation models with specific constraints include the \textit{1-2 model} \cite{grimmett_mot}, \textit{Eulerian percolation model} \cite{garet_mot}, and \textit{Random graphs with forbidden vertex degrees} \cite{grimmett_mot2}.

    In this paper, we propose and study a percolation model with restrictions on the difference between the volumes of the clusters adjacent to each site. We try to open the sites sequentially, the attempt will be successful if the volume difference between the two largest clusters adjacent to the site is greater than or equal to a constant $r$ (\textit{contrained}) or there is at most one cluster adjacent to it. The proposed model can be motivated by processes involving the spread and interaction of information or states within a networked system, where local dynamics are influenced by the relative sizes of interacting groups. For instance, such dynamics may occur in contexts like the diffusion of opinions, innovations, or behaviors, where an individual (site) adopts a state only if the influence of one group sufficiently outweighs that of another.
    
    This model will be called \textit{Constrained volume-difference site percolation model} and is formally defined below:

    \begin{figure}[t!]
        \centering
        \includegraphics[width=0.575\textwidth]{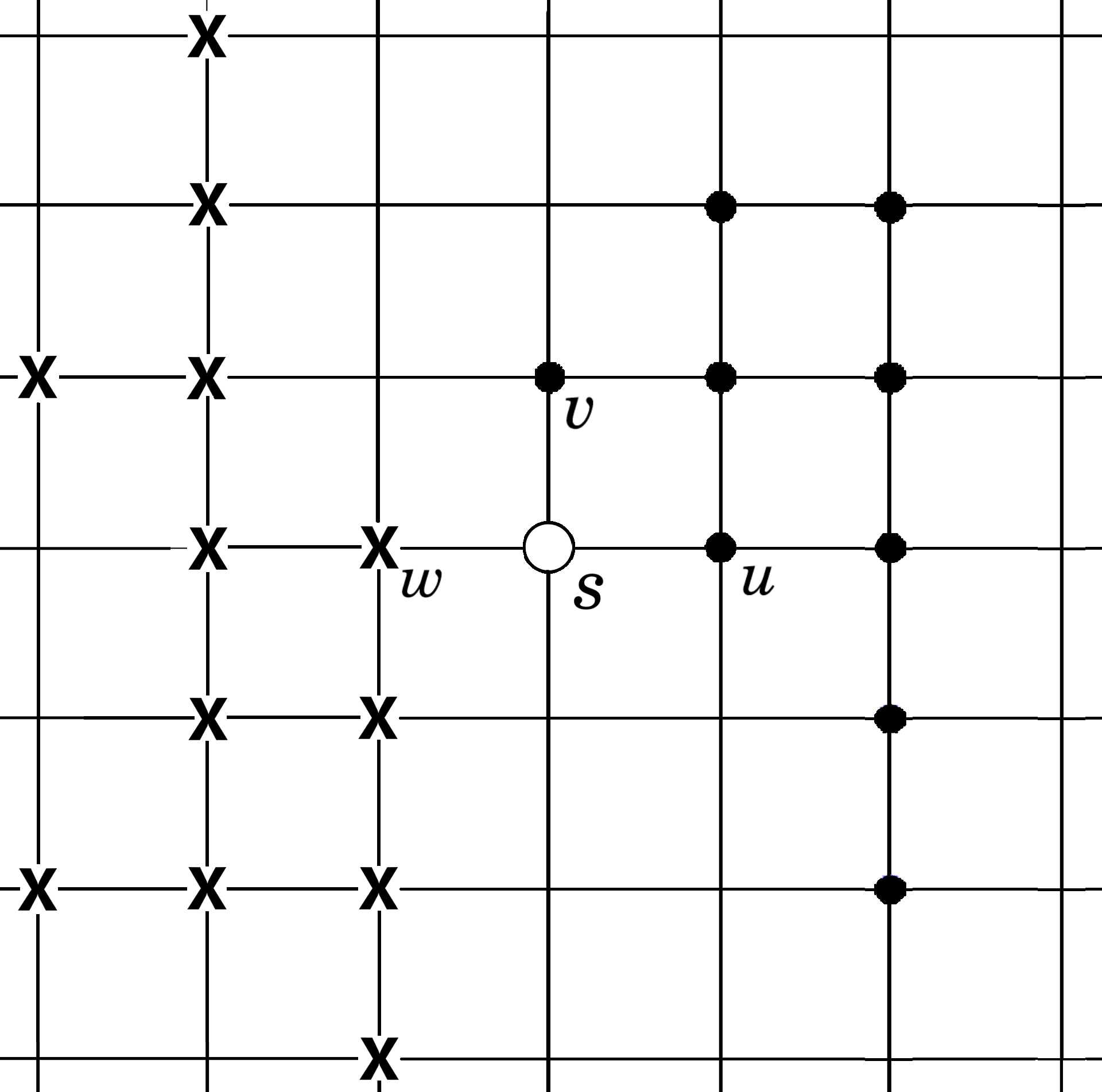}
        \caption{Part of the square lattice showing the attempt to open site $s$ at time $t=t_s$. All sites marked with the symbols $\bullet$ or $\times$ are open. Sites marked with identical symbols belong to the same cluster. Note that $M_{1}(t_s, s)=|\mathcal{C}_{w}|=12$ and $M_{2}(t_s, s)=|\mathcal{C}_{u}|=|\mathcal{C}_{v}|=9$. Therefore, since $|M_{1}(t_s, s) - M_{2}(t_s, s)| = 3$, the site $s$ will open if and only if the restriction $r$ chosen for the model satisfies $r \leq 3$.}
        \label{exemplo1}
    \end{figure}

    \begin{itemize}
        \item[(1)] Let $\mathcal{G} = (\mathcal{V}, \mathcal{E})$ be an infinite and connected graph with bounded degree. For each site $s \in \mathcal{V}$, assign a number $t_s$ randomly and uniformly chosen from the interval $[0,1]$.
        \item[(2)] Fix a non-negative integer $r$ and consider $t \in [0,1]$.
        \item[(3)] For each site $s$ and time $t$, denote by $M_{1}(t, s)$ and $M_{2}(t, s)$, $M_{1}(t, s) \geq M_{2}(t, s)$, the volumes of the two largest distinct clusters that contain some neighbor of site $s$. If a neighbor of site $s$ is closed at time $t$, then the volume of the cluster containing it is $0$.
        \item[(4)] Consider the continuous time process where at time $t = 0$ all sites are closed, and each site $s$ will open at time $t_s$ provided that $|M_{1}(t_s, s) - M_{2}(t_s, s)| \geq r$ or $M_{2}(t_s, s) = 0$. Note that if $M_{2}(t_s, s) = 0$, there will be at most one cluster connected to site $s$.
        \item[(5)] When $t = 1$, we have already attempted to open all the sites, and the model reaches the final state.
    \end{itemize}

   Note that for $r=0$, the condition $|M_1(t_s, s) - M_2(t_s, s)| \geq r$ is always satisfied. In this case, the model simplifies to OSPM, where sites open independently based on their assigned times $t_s$. For $r=1$, a site $s$ opens if the largest cluster adjacent to it exceeds the second-largest cluster by at least one site in volume, or if the second-largest cluster is empty ($M_2(t_s, s) = 0$). For $r=2$, this difference must be at least two sites, and so on for larger values of $r$.

   We study, via numerical simulations, this model on the square lattice $\mathbb{Z}^2$. In Figure \ref{exemplo1}, we show a part of this graph at the moment when we attempt to open the site $s$.

   If we replace the constraint $|M_{1}(t_s, s) - M_{2}(t_s, s)| \geq r$ by $|M_{1}(t_s, s) - M_{2}(t_s, s)| < r$, it seems that the model will not exhibit a phase transition for any value of $r$. This likely occurs because when there is a cluster significantly larger than the others, the constraint will prevent it from continuing to grow; we will only be able to open sites that have just this cluster connected to them. When the opposite constraint ($<r$) is considered, we simulated the model and found that percolation does not occur for $r \in \{1, \ldots, 9\}$.

    We denote by $\psi(t;r)$ the probability of percolation occurring at time $t$ when the constraint is equal to $r$. The critical time can be formally defined as:
    \begin{equation}
    t_{c}(r) = \sup\{ t \in [0,1]; \ \psi(t;r)=0 \}.
    \label{def_tc}
    \end{equation}
        	
    \noindent When $\psi(1;r)=0$, percolation does not occur for any value of $t$, and we define $t_{c}(r)~=~\infty$. If $r=0$, this model reduces to OSPM with parameter $p=t$. Therefore, $t_c(0)=p_c(\mathbb{Z}^2)=0.592746050786(3)$ \cite{mertens2} (the ordinary site percolation threshold).
    
    When the graph is the one-dimensional lattice $\mathbb{Z}$, we have that $t_c(0)=p_{c}(\mathbb{Z})=1$ and $t_c(r)=\infty$ for any $r \geq 1$. In fact, let $(t_{s})_{s \in \mathbb{Z}}$ be the sequence of random numbers where $t_s$ is the time we attempt to open site $s$. For any site $k$, we have that the probability of the event
    $$A_{k}:=\{t_k > t_{k-1} < t_{k-2}\} \cap \{t_k > t_{k+1} < t_{k+2}\} \cap \{t_k < t_{k-2}\} \cap \{t_k < t_{k+2}\}$$ 
    \noindent is positive, $\mathbb{P}(A_k)=\delta>0$, and its occurrence implies that when we attempt to open site $k$, $t=t_{k}$, we will have $|\mathcal{C}_{k-1}|=|\mathcal{C}_{k+1}|=1$. Therefore, $|M_{1}(t_{k}, k) - M_{2}(t_{k}, k)| = 0$ and site $k$ will not open. Since the event $A_k$ depends only on the five times $t_{k-2}, t_{k-1}, t_k, t_{k+1}, t_{k+2}$, the events $(A_{5k})_{k \geq 0}$ are independent. Furthermore, since $\sum_{k \geq 0} \mathbb{P}(A_{5k}) = \sum_{k \geq 0} \delta = \infty$, the second Borel-Cantelli Lemma (see \cite{shiryaev}) ensures that infinitely many events $A_{5k}$ occur with probability $1$, thereby preventing percolation for any $t \in [0,1]$.

    The higher the value of $r$, it is intuitive to expect that we will have fewer open sites. Therefore, we conjecture that if percolation occurs for $r = r_1$ and $r = r_2$, with $r_1 < r_2$, then $t_{c}(r_1) \leq t_{c}(r_2)$. In other words, $t_{c}(r)$ is non-decreasing with $r$. Furthermore, this argument suggests the existence of a critical value $r_c$, such that for $r > r_c$, percolation does not occur and, thus, the model does not exhibit a phase transition ($t_c(r) = \infty$).

    Through numerical simulations of this model, we estimate the critical time for several values of $r$, finding evidence that $t_{c}(r)$ is indeed non-decreasing with $r$ and that $r_c = 5$. We also found that the correlation length exponent $\nu$ is equal to that of the ordinary percolation model. Furthermore, we determine other properties of the model when $t=1$ and $1 \leq r \leq 9$: the average density of open sites ($\overline{\rho}(r)$), the average number of distinct cluster volumes ($\overline{n}(r)$), and the average volume of the largest cluster ($\overline{M}(r)$).
                
    In the Constrained volume-difference site percolation model, when $r > 0$, the probability of a site being open or closed depends on the status of other sites (long-range dependence), unlike the OSPM and OBPM. Consequently, algorithms like the Leath Algorithm \cite{leath} or the Invasion Percolation Algorithm \cite{chandler, wilkinson2}, which assume the independence of each site being open or closed, cannot be used.  These algorithms generate only a single cluster rather than an entire lattice configuration and are often employed to estimate the critical point of percolation models due to their low computational time and space complexities.
    
    The remainder of this paper is organized as follows. Section 2 describes our simulation procedure and Section 3 discusses the results obtained. Conclusions are summarized in Section 4.

    \section{Numerical Procedure} \label{num_procedure}
	
    To estimate $t_{c}(r)$, we consider periodic square lattices of lengths $L=2^j$, where $j$ ranged from $8$ to $12$, to be the set of sites of the graph. For each constraint $r \in \{1,...,9\}$, we use the Newman-Ziff algorithm \cite{ziff01} to estimate the critical time $t_c(r)$. The percolation criterion assumed is that the infinite cluster emerges when there exists a cluster that wraps around either the horizontal or vertical directions. The number of simulations ranged from $10^5$ ($L=4096$) to $5 \times 10^6$ ($L=256$). The parameters $\overline{\rho}(r)$, $\overline{n}(r)$ and $\overline{M}(r)$ were estimated considering only $L=2048$ and $5 \times 10^ 4$ simulations.

    We denote by $\psi_{L}(t;r)$ the probability of percolation in the square of length $L$. For each configuration $ (t_s)_{s \in \mathbb{Z}^2} $, we define $\mathcal{O}_t$ as the number of sites that we will try to open until time $t$, that is  
        \begin{equation}
        \mathcal{O}_t = \# \{s \in \mathbb{Z}^2; t_s \leq t \}
        \end{equation}
    \noindent where the $\#$ symbol denotes the number of elements in the set. Note that the probability of a site belongs to $\mathcal{O}_t$ does not depend on the state of any other site (the same does not occur when we consider the set formed only by sites opened up to time $t$). This independence allows us to calculate explicitly the probability that $\mathcal{O}_t=i$, for all $ i $, and we can write, using the \textit{Partition Theorem} (see \cite{grimmett_prob}),  
                	
    \begin{equation}
    \psi_{L}(t;r)= \sum_{i=0}^{N} \mathbb{P} (\mathcal{O}_t=i) \cdot \overline{Q_i} = \sum_{i=0}^{N} 
    \left(\begin{array}{c}
    N \\
    i 
    \end{array}\right) \cdot t^{i} \cdot (1-t)^{(N-i)} \cdot \overline{Q_i}
    \label{psi_L}
    \end{equation}

    \noindent where $N=L^{2}$ (number of sites in the graph) and $\overline{Q_i}$ is the probability of percolation, at time $t$, conditioned that $\mathcal{O}_t=i$. 

    To determine $\overline{Q}_{i}$, we create a list with a permutation of all sites and attempt to open them in the obtained order. If a site satisfies the constraint, we open it and check if percolation has occurred. If percolation occurs when the $i$-th site is opened, then we define $Q_j = 1$ for $j \geq i$. Otherwise, $Q_i = 0$. The estimate of $\overline{Q}_{i}$, for each $i$, is obtained as the mean of many samples of $Q_i$.

    For each $r$, the critical time, $t_{c}(r)$, is obtained by the finite-size scaling (FSS)
    \begin{equation}
    |\overline{t}_{L}(r)-t_{c}(r)| \sim L^{-\frac{1}{\nu(r)}},
    \label{tc}
    \end{equation}
            	
    \noindent where $\overline{t}_{L}(r) = \int_{0}^{1} t \cdot \frac{d\psi_{L}(t; r)}{dt} \, dt$ is the average concentration when percolation occurs for the first time, and $\nu(r)$ is the \textit{correlation length exponent} \cite{stauffer94}. The critical exponent $\nu(r)$ was estimated from the scaling relation  
\begin{equation}
    \mathrm{Max}\left(\psi'_{L}(t; r)\right) \sim L^{\frac{1}{\nu(r)}},
    \label{nu_fss}
\end{equation}  
\noindent where the left-hand side represents the maximum derivative of $\psi_{L}(t; r)$ with respect to $t$. The uncertainties in $t_{c}(r)$ and $\nu(r)$ were determined from the standard deviation of regression residuals.

    In the estimation of $\overline{\rho}(r)$, $\overline{n}(r)$, and $\overline{M}(r)$, we analyzed the final configuration in each of the $5 \times 10^4$ simulations. The uncertainties were determined using the standard deviation across all simulations.
    
    This computational task consumed approximately $1.5$ months of CPU time across 20 cores operating at a clock speed of 2.7 GHz.

\section{Results}

    The values of $\nu(r)$, presented in Table \ref{table_tc}, are consistent with the correlation length exponent for ordinary percolation on the square lattice, $\nu = \frac{4}{3}$ \cite{stauffer94}, providing evidence that both models belong to the same universality class. 

    \begin{table}[h!]
        \renewcommand{\arraystretch}{1.2} 
        \setlength{\extrarowheight}{0.1cm} 
        \setlength{\tabcolsep}{3.1pt} 
        \centering
        \begin{tabular}{c|c|c|c|c|c|c|}
            \cline{2-7}  
            & \multicolumn{1}{c|}{\textbf{$r=1$}} & \multicolumn{1}{c|}{\textbf{$r=2$}} & \multicolumn{1}{c|}{\textbf{$r=3$}} & \multicolumn{1}{c|}{\textbf{$r=4$}} &  \multicolumn{1}{c|}{\textbf{$r=5$}} & \multicolumn{1}{c|}{\textbf{$6 \leq r \leq 9$}} \\ \cline{2-6} \hline
            \multicolumn{1}{|c|}{\textbf{$t_c(r)$}}    & 0.633306(4)        & 0.701872(8)        & 0.77913(2)         & 0.86437(3)        & 0.95839(7)  & $\infty$ 	 \\ \hline
            \multicolumn{1}{|c|}{\textbf{$\nu(r)$}}    & 1.34(1)        & 1.34(1)        & 1.34(1)         & 1.34(1)        & 1.33(1)  & $-$ 	 \\ \hline
        \end{tabular}
        \caption{Estimated critical times $t_{c}(r)$, for $1\leq r \leq 9$, and $\nu(r)$, for $1\leq r \leq 5$.}
        \label{table_tc}
    \end{table}
    
\begin{figure}[t!]
        \centering
        \includegraphics[width=0.65\textwidth, height=9.0cm]{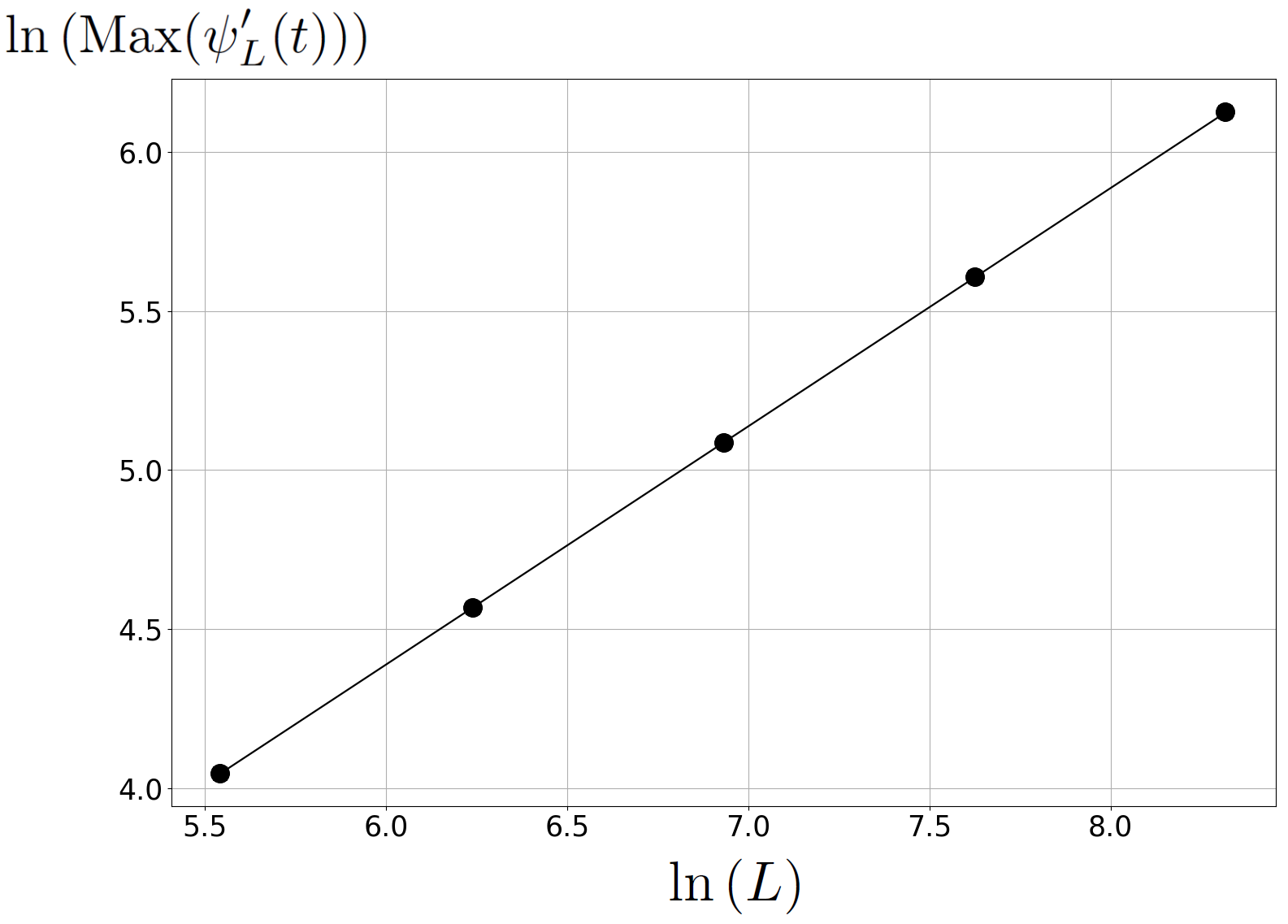}
        \caption{Graph obtained through a log-log plot using the FSS relation (\ref{nu_fss}) for the case $r=1$. The slope of the linear regression corresponds to the estimated value of $\frac{1}{\nu(r)}$.}
        \label{fig_nu}
    \end{figure}

    The graph obtained through the relation (\ref{nu_fss}) to estimate $\nu(r)$ for $r=1$ is shown in Fig.~\ref{fig_nu}. Graphs for $2 \leq r \leq 5$ are omitted as they are qualitatively similar to the case $r=1$.

\begin{figure}[t!]
        \centering
        \includegraphics[width=0.65\textwidth, height=8.7cm]{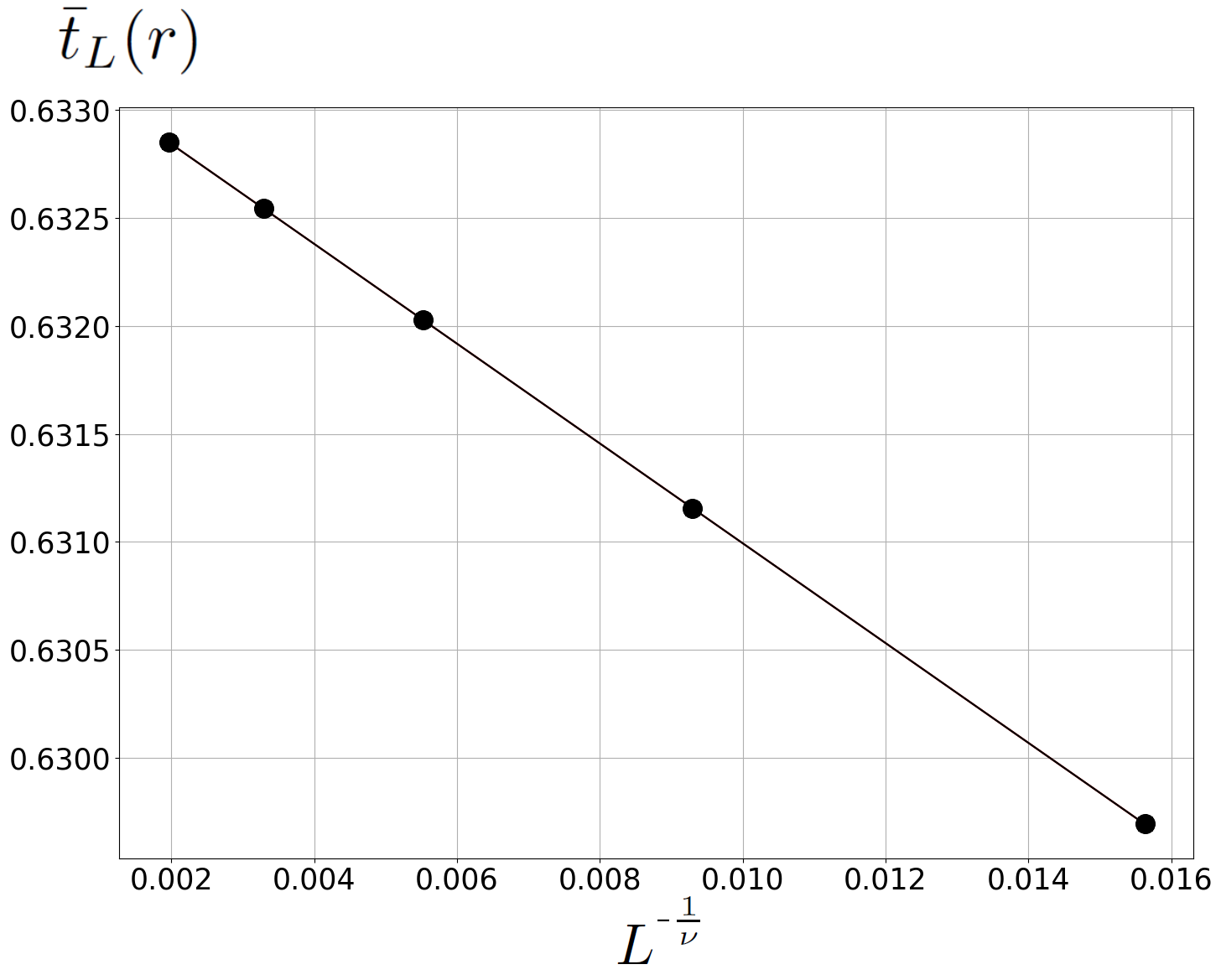}
        \caption{Graph obtained through the FSS relation (\ref{tc}) to estimate the critical time $t_c(r)$ for the case $r=1$. The linear coefficient of the line obtained from the linear fit is the estimate for $t_c(1)$.}
        \label{fig_tc}
    \end{figure}

\begin{figure}[t!]
        \centering
        \includegraphics[width=0.63\textwidth, height=9.3cm]{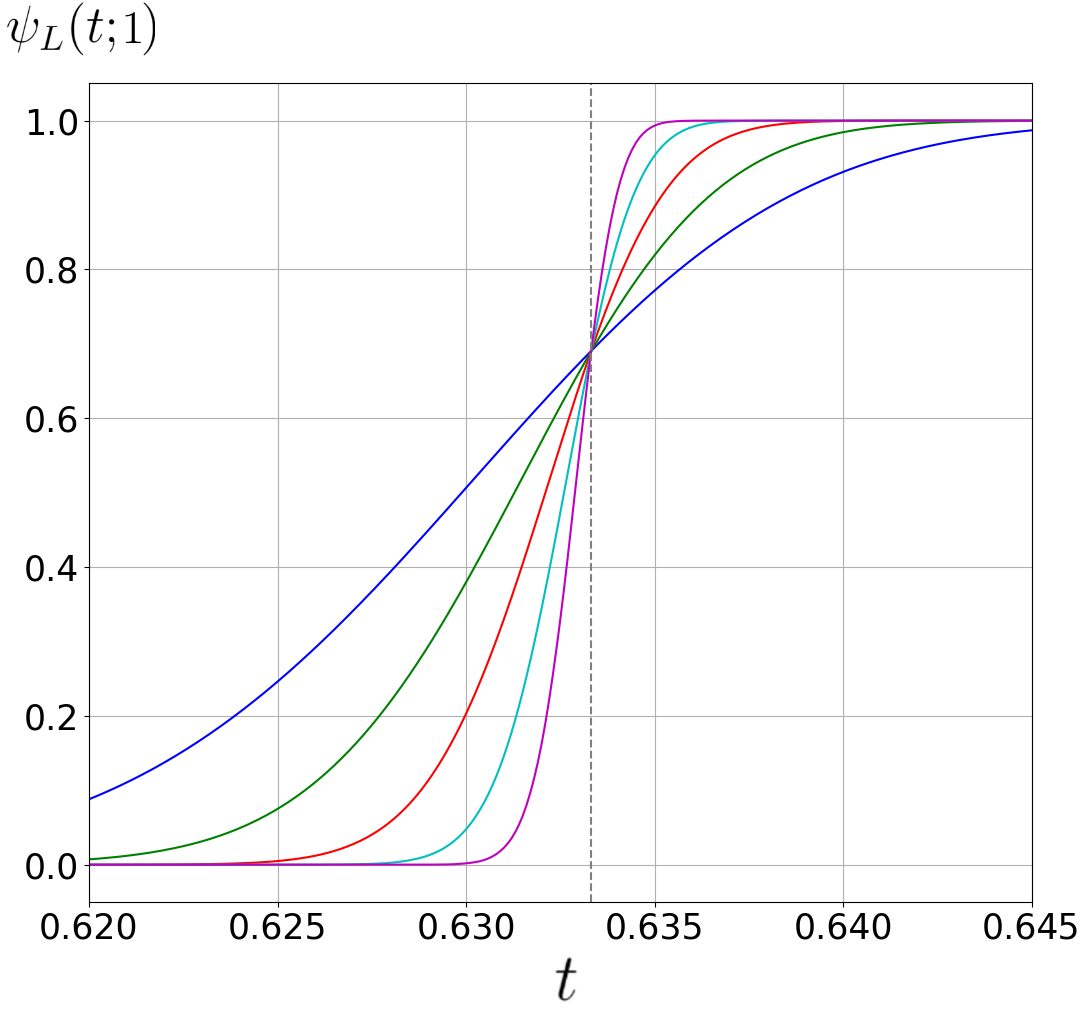}
        \caption{Graphs of the $\psi_{L}(t;r)$ functions for $r=1$ and all values of $L$. The dashed line indicates the estimated critical time $t_c(r) = 0.633306(4)$. The slope of the graph near the critical time increases as $L$ grows.}
        \label{fig_psi}
    \end{figure}

    The critical times $t_c(r)$ were estimated using the FSS relation (\ref{tc}), and the graph for $r=1$ is shown in Fig. \ref{fig_tc}. In Fig. \ref{fig_psi}, we present the behavior of $\psi_{L}(t;1)$ for all $L$, where we observe that the curves intersect near the critical value $t_c(1) = 0.633306(4)$. 

    We assume $\nu(r) = \frac{4}{3}$ for all $r$ and the estimated values of $t_c(r)$ are presented in Table 1. For $r>5$, the model does not exhibit a phase transition ($t_c(r) = \infty$), hence $r_c = 5$. This absence of phase transition highlights the significant impact of the constraint parameter $r$ on the percolation behavior of the system. 
    
    We also observed that $t_{c}(r)$ is non-decreasing with $r$ and exhibits an almost linear growth for $1 \leq r \leq 5$, as $\frac{t_c(r+1)}{t_c(r)} \approx 1.1$ for $1 \leq r \leq 4$.

    The results obtained for $t_c(r)$ indicate that the system becomes more resistant to percolation as $r$ increases, requiring longer times for it to occur when $r \leq r_c$. This behavior is consistent with the expectation that stricter constraints hinder the formation of a giant component in the system.
    
    To further investigate the effects of the constraint parameter $r$, we analyzed the system at $t=1$ and estimated the average density of open sites $\overline{\rho}(r)$, the average proportion of nodes in the largest cluster $\overline{M}(r)$, and the average number of distinct cluster volumes $\overline{n}(r)$. The results are summarized in Table 2.

    \begin{table}[t!]
        \renewcommand{\arraystretch}{1.3}
        \setlength{\extrarowheight}{0.1cm}
        \setlength{\tabcolsep}{6pt}
        \centering
        \begin{tabular}{c|c|c|c|}
        \cline{2-4}
        & \multicolumn{1}{c|}{$\overline{\rho}(r)$} & \multicolumn{1}{c|}{$\overline{M}(r)$ (\mbox{\%})} & \multicolumn{1}{c|}{$\overline{n}(r)$}  \\
        \cline{2-4}
        \hline
        \multicolumn{1}{|c|}{$ r=1 $} & 0.9672(2) & 96.69(2) & 4(1) \\
        \multicolumn{1}{|c|}{$ r=2 $} & 0.9107(3) & 90.81(3) & 11(2) \\
        \multicolumn{1}{|c|}{$ r=3 $} & 0.8463(4) & 83.29(6) & 26(2) \\
        \multicolumn{1}{|c|}{$ r=4 $} & 0.7770(6) & 72.88(11) & 70(4) \\
        \multicolumn{1}{|c|}{$ r=5 $} & 0.7062(8) & 56.26(30) & 252(11) \\
        \multicolumn{1}{|c|}{$ r=6 $} & 0.6437(8) & 0.54(14) & 1039(17) \\
        \multicolumn{1}{|c|}{$ r=7 $} & 0.5964(7) & 0.09(2) & 867(12) \\
        \multicolumn{1}{|c|}{$ r=8 $} & 0.5652(6) & 0.04(1) & 585(10) \\
        \multicolumn{1}{|c|}{$ r=9 $} & 0.5467(5) & 0.02(1) & 428(8) \\
        \hline
        \end{tabular}
        \caption{Values obtained for $\overline{\rho}(r)$, $\overline{M}(r)$, and $\overline{n}(r)$ considering $5 \times 10^4$ simulations of the model for $L=2048$.}
        \label{t1}
    \end{table}

    As $r$ increases, it becomes less likely for each site to open and, therefore, $\overline{\rho}(r)$ decreases. It is worth noting that our model can be interpreted as a \textit{Random Sequential Adsorption} model (RSA). In such models, particles are sequentially added to a surface according to specific rules, and the probability of a particle adsorbing onto the surface depends on the spatial arrangement of previously adsorbed particles \cite{evans1993}. An important feature of RSA is the \textit{jamming limit}, which corresponds to the average density of open sites on the graph at the end of the deposition process. In the Constrained volume-difference site percolation model, this limit is precisely given by the value of $\overline{\rho}(r)$.

    The percentage of sites in the largest cluster provides insights into the formation of dominant clusters. For $r=1$, $\overline{M}(r)$ comprises $96.69(2)\%$ of all the sites in the graph. As $ r $ increases, $\overline{M}(r)$ decreases drastically, with the largest cluster comprising only $0.02(1)\%$ of the sites for $r=9$. This trend underscores the fragmentation effect caused by higher constraints, leading to smaller and more dispersed clusters. The most significant decrease in $\overline{M}(r)$ occurs between the values $r=5$ and $r=6$. This behavior is explained by the fact that percolation does not occur for $r>r_c=5$.
    
    The number of distinct cluster volumes $\overline{n}(r)$ reveals the diversity of cluster formations within the network. For $ r = 1 $, $\overline{n}(r)$ is relatively low, $3.8(6)$, indicating that a few distinct cluster volumes dominate the network. This occurs because the largest cluster dominates nearly the entire graph, leaving only smaller clusters in addition to it.
    
    For $1 \leq r \leq 6$, $\overline{n}(r)$ increases and shows the largest jump between $r=5$ and $r=6$, as the absence of percolation for $r > r_c$ allows for a greater diversity of cluster volumes. For $r > 6$, the constraint leads to a decrease in the volume of the largest cluster, $\overline{M}(6)=0.54$ to $\overline{M}(7)=0.09$, resulting in only smaller clusters in the graph and thus a decrease in $\overline{n}(r)$.

 \section{Conclusion}

In conclusion, we have introduced and analyzed the Constrained volume-difference site percolation model on the square lattice. In this model, starting with all sites closed, we attempt to open them sequentially. The attempt will be successful if the volume difference between the two largest clusters adjacent to the site is greater than or equal to a constant $r$ or if there is at most one cluster adjacent to the site. Our numerical simulations revealed that the critical time $t_c(r)$ increases with the constraint parameter $r$, for $1 \leq r \leq 5$, and identified a critical value $r_c = 5$, beyond which percolation does not occur. We also found that the correlation length exponent $\nu$ of this model is equal to that of the ordinary percolation model. 

Additionally, the average density of open sites $\overline{\rho}(r)$ and the average proportion of nodes in the largest cluster $\overline{M}(r)$ decrease with increasing $r$, while the average number of distinct cluster volumes $\overline{n}(r)$ initially increases up to $r=6$ and then decreases for higher values of $r$. These results confirm that stricter constraints hinder the formation of large clusters in the system.

These findings contribute to an understanding of how local constraints influence global percolation properties. The behavior observed in this model, especially the absence of a phase transition for high constraint values, underscores the importance of constraint dynamics in network formation. Future work may explore the study of this model in higher-dimensional hypercubic lattices and other graph topologies.

    \section*{Acknowledgements}

    We thank Nuno A. M. Araujo and Diogo C. dos Santos for their helpful suggestions and comments that improved the manuscript. Charles S. do Amaral was partially supported by FAPEMIG (Processo APQ-00272-24).

\end{document}